\documentclass[12pt]{article}%
\usepackage{amsfonts, amsmath, amscd}
\usepackage[psamsfonts]{amssymb}

\newcounter{df}%[section]                      %
\newcounter{pro}%[section]                      %
\newcounter{rem}%[section]                     %
\newenvironment{rem}{\par%                    %
\refstepcounter{rem}%                         % оформление замечаний
{\bf Remark \arabic{rem}.} }{}%  \arabic{section}.       %
\newcounter{exa}%[section]                     %
\newcounter{teo}%[section]                      %
\newenvironment{teo}{\par%                     %
\refstepcounter{teo}%                          % оформление теорем
{\bf Theorem \arabic{teo}.} \it }{}% \arabic{section}.       %
\newcounter{cor}%[section]                     %
\newcounter{st}%[section]                      %
\newcounter{lem}%[section]                     %
\newenvironment{lem}{\par%                    %
\refstepcounter{lem}%                         % оформление лемм
{\bf Lemma \arabic{lem}.} \it }{}%  \arabic{section}.       %
\makeatletter                                                   %
\renewcommand{\section}{\@startsection{section}{1}%             % оформление
{\parindent}{3.5ex plus 1ex minus 0.2ex}{2.3ex plus 0.2ex}{\bf}}% заголовков
\makeatother                                                    %
\begin{document}%

\author{S.S. Gabriyelyan\footnote{The author was partially supported
 by Israel Ministry of Immigrant Absorption}}
\title{Characterization of almost maximally almost-periodic groups}
\date{}

\makeatletter
\renewcommand{\@makefnmark}{}
\renewcommand{\@makefntext}[1]{\parindent=1em #1}
\makeatother

\maketitle

\vspace{-\baselineskip}

\begin{abstract}
Let $G$ be an  abelian group. We prove that a group $G$ admits a Hausdorff group topology $\tau$ such that the von Neumann radical $\mathbf{n}(G, \tau)$ of $(G, \tau)$ is non-trivial and finite iff $G$ has a non-trivial finite subgroup. If $G$ is a topological group, then $\mathbf{n} (\mathbf{n} (G)) \not= \mathbf{n} (G)$ if and only if $\mathbf{n} (G)$ is not dually embedded. In particular, $\mathbf{n} (\mathbf{n} (\mathbb{Z},\tau)) = \mathbf{n} (\mathbb{Z},\tau)$ for any Hausdorff group topology $\tau$ on $\mathbb{Z}$.
\end{abstract}

We shall write our abelian groups additively. For a topological group $X$, $\widehat{X}$ denotes the group of all continuous characters on $X$. We denote its dual group by $X^{\wedge}$, i.e. the group $\widehat{X}$ endowed with the compact-open topology.  Denote by $\mathbf{n}(X) = \cap_{\chi\in \widehat{X}} {\rm ker} \chi$ the von Neumann radical of $X$. If $H$ is a subgroup, we denote by $H^{\perp}$ its annihilator.  If $A$ be a subset of a group $X$, $\langle A\rangle$ denotes the subgroup generated by $A$.

Let $X$ be a topological group and $\mathbf{u} =\{ u_n \}$ a sequence of elements of $\widehat{X}$. We denote by $s_{\mathbf{u}} (X)$ the set of all $x\in X$ such that $(u_n , x)\to 1$. Let $G$ be a subgroup of $X$. If $G=s_{\mathbf{u}} (X)$ we say that $\mathbf{u}$ {\it characterizes} $G$ and that $G$ is {\it characterized} (by $\mathbf{u}$).

Following E.G.Zelenyuk and I.V.Protasov \cite{ZP1}, \cite{ZP2}, we say that a sequence $\mathbf{u} =\{ u_n \}$ in a group $G$ is a $T$-{\it sequence} if there is a Hausdorff group topology on $G$ for which $u_n $ converges to zero. The group $G$ equipped with the finest group topology with this property is denoted by $(G, \mathbf{u})$.

Group topologies on $\mathbb{Z} (p^{\infty})$ with $\mathbf{n} (\mathbb{Z} (p^{\infty})) = \mathbb{Z} (p)$ were considered in corollary 4.9 \cite{DMT}. Although no explicit construction of such topology was given there, it was conjectured by D.~Dikranjan that such topology can be found by means of an appropriate $T$-sequence on $\mathbb{Z} (p^{\infty})$. This conjecture was successfully proved by G.~Luk\'{a}cs  \cite{Luk} for every prime $p\not= 2$. He called a Hausdorff topological group $G$ {\it almost maximally almost-periodic} if $\mathbf{n} (G)$ is non-trivial and finite and raised the problem of their description. He proved that infinite direct sums and the Pr\"{u}fer group $\mathbb{Z} (p^{\infty}),$ for every prime $p\not= 2$, are almost maximally almost-periodic. A.P.Nguyen  \cite{Ngu} generalized these results and proved that any Pr\"{u}fer groups $\mathbb{Z} (p^{\infty})$ and a wide class of torsion groups admit a (Hausdorff) almost maximally almost-periodic group topology. Using theorem 4 \cite{Ga1}, we give a general characterization of almost maximally almost-periodic groups.
\begin{teo} \label{t1}
{\it Let $G$ be an infinite group. Then the following statements are equivalent.
\begin{enumerate}
\item $G$ admits a $T$-sequence $\mathbf{u}$ such that $(G, \mathbf{u})$ is almost maximally almost-periodic.
\item $G$ has a non-trivial finite subgroup.
\end{enumerate} }
\end{teo}

Evidently that $\mathbf{n} (\mathbf{n} (G)) \not= \mathbf{n} (G)$ if $\mathbf{n} (G)$ is non-trivial and finite. This observation leaded G. Luk\'{a}cs \cite{Luk} to the problem of description of abelian groups $G$ which admit a $T$-sequence $\mathbf{u}$ such that $\mathbf{n} (\mathbf{n} (G,\mathbf{u}))$ is strictly contained in $\mathbf{n} (G,\mathbf{u})$. We can generalize this question in the following way.
\begin{enumerate}
\item[] {\bf Problem.} Which abelian groups $G$ admit a Hausdorff group topology $\tau$ such that $\mathbf{n} (\mathbf{n} (G,\tau)) \not= \mathbf{n} (G,\tau)$?
\end{enumerate}

By theorem \ref{t1}, if $G$ is a not torsion free group, the answer is positive. But for torsion free groups an answer is negative in general. We give an answer to this problem in the next theorem.

\begin{teo} \label{t2}
{\it Let $G$ be a topological group. Then $\mathbf{n} (\mathbf{n} (G)) \not= \mathbf{n} (G)$ if and only if $\mathbf{n} (G)$ is not dually embedded. In particular, $\mathbf{n} (\mathbf{n} (\mathbb{Z},\tau)) = \mathbf{n} (\mathbb{Z},\tau)$ for any Hausdorff group topology $\tau$ on $\mathbb{Z}$.}
\end{teo}

\begin{center}
{\bf\large The proofs}
\end{center}

The following lemma plays an important role for the proofs of theorems \ref{t1} and \ref{t2}.
\begin{lem} \label{l3}
{\rm \cite{Ga1}} {\it Let $H$ be a dually closed and dually embedded subgroup of a topological group $G$. Then $\mathbf{n} (H) =\mathbf{n} (G)$. }
\end{lem}

{\bf Proof of theorem \ref{t1}.} It is clear that $1\Rightarrow 2$.

Let us prove that $2\Rightarrow 1$. Denote by $H$ any countable subgroup which contains an element of finite order. If $H$ admits an almost maximally almost-periodic group topology, then we can consider the topology on $G$ in which $H$ is open. Since any open subgroup is dually closed and dually embedded (lemma 3.3 \cite{Nob}), by lemma \ref{l3},
$G$  admits an almost maximally almost-periodic group topology too. Thus we can and will pass to countable subgroups. In particular, we will suppose that $G$ is countable.
There exist three possibilities.

1) $G$ contains an element $g$ of infinite order and an element $e$ of order $p$ for some prime $p$. Thus $G$ contains the subgroup $\langle e , g\rangle $ which is isomorphic to $\mathbb{Z} (p) \oplus \mathbb{Z}$. By lemma \ref{l3},
we can assume that $G=\mathbb{Z} (p) \oplus \mathbb{Z}$. Then $G^{\wedge} =\mathbb{Z} (p) \oplus \mathbb{T}$. Now we define the following sequence
$$
d_{2n-1} = e + p^{n^3 - n^2} +\dots + p^{n^3 - 2n} +  p^{n^3 - n} + p^{n^3}, \quad d_{2n} =  p^{n}, \quad n\geq 1.
$$
Let us prove that $\mathbf{d} =\{ d_n\}$ is a $T$-sequence.

For our convenience, we put $f_n =p^{n^3 - n^2} +\dots + p^{n^3 - 2n} +  p^{n^3 - n} + p^{n^3}$, then $f_n < 2 p^{n^3} \leq p^{n^3 +1}$. For any $0<r_1 <r_2 <\dots <r_v$ and integers $l_1, l_2,\dots , l_v$ such that $\sum_{i=1}^v |l_i | \leq k+1$, we have
\begin{equation} \label{3}
| l_1 f_{r_1 } + l_2 f_{r_2 } +\dots + l_v f_{r_v } | < (k+1) f_{r_v} \leq (k+1) p^{r_v^3 +1} ,
\end{equation}

Let $g= ae+b \not= 0, 0\leq a <p, b\in \mathbb{Z},$ and $k\geq 0$. Set $t=|b|+p(k+1)$ and $m=10t$. We shall prove that $g\not\in A(k,m)$.

a) Let $\sigma \in A(k,m)$ have the following form
$$
\sigma =l_1 d_{2r_1 } + l_2 d_{2r_2 } +\dots + l_s d_{2r_s} = l_1 p^{r_1}  +\dots + l_s p^{r_s} = p^{r_1} \cdot\sigma' ,
$$
where  $m\leq 2r_1 <2r_2 <\dots <2r_s$  and  $\sigma' \in \mathbb{Z}$.
If $\sigma' =0$, then $\sigma \not= g$. If $\sigma' \not= 0$, then $|\sigma |\geq p^{r_1} \geq p^{5|b|} >|b|$, and $\sigma \not= g.$

b) Let $\sigma \in A(k,m)$ have the following form
$$
\sigma =l_1 d_{2r_1 -1} + l_2 d_{2r_2 -1} +\dots + l_s d_{2r_s -1},
$$
where $m< 2r_1 -1 < 2r_2 -1 < \dots < 2r_s -1 $ and integers $l_1, l_2,\dots , l_s$ be such that $l_s \not= 0$ and  $\sum_{i=1}^s |l_i | \leq k+1$. Since $n^3 < (n+1)^3 -(n+1)^2$ and $r_s >5p(k+1)$, by (\ref{3}), we have
$$
|\sigma -(l_1 +\dots +l_s) e | = |l_1 f_{r_1 } + \dots + l_{s-1} f_{r_{s-1} } + l_s f_{r_s } | >
f_{r_s } - (k+1) p^{r_{s-1}^3 +1} =
$$
$$
p^{r_{s}^3 } +\dots + \left( p^{r_{s}^3 -r_s (r_s -1)} + p^{r_{s}^3 -r_s^2} - (k+1) p^{r_{s-1}^3 +1} \right) > \left(\mbox{since } (k+1)p < p^{r_{s}-1} \right) > p^{r_{s}^3 } >|b|
$$
and $\sigma \not= g$.

c) Let $\sigma \in A(k,m)$ have the following form
$$
\sigma = l_1 d_{2r_1 -1} + l_2 d_{2r_2 -1} +\dots + l_s d_{2r_s -1} +l_{s+1} d_{2r_{s+1}} +l_{s+2} d_{2r_{s+2}} +\dots + l_{h} d_{2r_{h}},
$$
where $0<s<h$ and
$$
\begin{array}{lr}
  m< 2r_1 -1 < 2r_2 -1 < \dots < 2r_s -1 , &   \\
  m\leq 2r_{s+1} < 2r_{s+2} < \dots < 2r_{h}, & l_i \in \mathbb{Z} \setminus \{ 0\} , \; \sum_{i=1}^h |l_i | \leq k+1.
\end{array}
$$
Since the number of summands in $d_{2r_s -1}$ is $r_s +2> 5p(k+1)+2$ and $h-s <k+1$, there exists $r_s - 2> i_0 >2$ such that for every $1\leq w \leq h-s$ we have
$$
\mbox{either } r_{s+w} < r_s^3 - (i_0 +2)r_s  \mbox{ or }  r_{s+w} > r_s^3 - (i_0 -1)r_s .
$$
The set of all $w$ such that $r_{s+w} < r_s^3 - (i_0 +2)r_s$ we denote by $B$ (it may be empty or has the form $\{ 1,\dots ,\delta \}$ for some $1\leq \delta \leq h-s$). Set $D =\{ 1,\dots ,h-s \} \setminus B$. Thus
$$
\begin{array}{ll}
  \sigma =  & (l_1 +\dots +l_s) e + l_1 f_{r_1 } + \dots + l_{s-1} f_{r_{s-1} }  +\\
   & \sum_{w\in B} l_{s+w} d_{2r_{s+w}} + l_s p^{r_s^3 - r_s^2} +\dots + l_s p^{r_s^3 - (i_0 +2)r_s} + \\
   &  l_s p^{r_s^3 - (i_0 +1)r_s} + l_s  p^{r_s^3 - i_0 r_s} +  \\
   &  l_s p^{r_s^3 - (i_0 -1)r_s} +\dots + l_s  p^{r_s^3} + \sum_{w\in D} l_{s+w} d_{2r_{s+w}}.
\end{array}
$$
We can estimate the expression in row 2, which we denote by $A_2$, as follows
\begin{equation} \label{4}
|A_2 | < \sum_{w\in B} |l_{s+w} | p^{r_s^3 - (i_0 +2)r_s}+ |l_s | 2 p^{r_s^3 - (i_0 +2)r_s} <3(k+1) p^{r_s^3 - (i_0 +2)r_s} <p^{r_s^3 - (i_0 +1)r_s},
\end{equation}
since $3(k+1) <r_s < p^{r_s} -1$. For the expression in row 4, which we denote by $A_4$, we have
\begin{equation} \label{5}
A_4 = l_s p^{r_s^3 - (i_0 -1)r_s} +\dots + l_s  p^{r_s^3} + \sum_{w\in D} l_{s+w} d_{2r_{s+w}} =  p^{r_s^3 - (i_0 -1)r_s} \cdot \sigma'' ,
\end{equation}
where  $\sigma'' \in \mathbb{Z}$. By (\ref{3})- (\ref{5}), we have: if $\sigma'' \not= 0$, then
$$
\begin{array}{ll}
  |\sigma -(l_1 +\dots +l_h) e |= & |l_1 f_{r_1 } + \dots + l_{s-1} f_{r_{s-1} } + A_2 + l_s p^{r_s^3 - (i_0 +1)r_s} + l_s  p^{r_s^3 - i_0 r_s} +  A_4 | > \\
   & p^{r_s^3 - (i_0 -1)r_s } -(k+1)p^{r_{s-1}^3 +1} -2(k+1) p^{r_s^3 - i_0 r_s} > \\
   & p^{r_s^3 - (i_0 -1)r_s } -3(k+1) p^{r_s^3 - i_0 r_s} > p^{r_s^3 - i_0 r_s }>  p^{r_s^2} > p^{|b|} >|b|
\end{array}
$$
and $\sigma \not= g$; if $\sigma'' = 0$, then, by (\ref{3}) and (\ref{4}),
$$
\begin{array}{ll}
  |\sigma -(l_1 +\dots +l_h) e |= & |l_1 f_{r_1 } + \dots + l_{s-1} f_{r_{s-1} } +A_2  + l_s p^{r_s^3 - (i_0 +1)r_s} + l_s  p^{r_s^3 - i_0 r_s} | > \\
   & p^{r_s^3 - i_0 r_s} - (k+1)p^{r_{s-1}^3 +1} -(k+2) p^{r_s^3 - (i_0 +1)r_s} > \\
   & p^{r_s^3 - i_0 r_s} -3(k+1) p^{r_s^3 - (i_0 +1)r_s } > p^{r_s^3 - (i_0 +1)r_s } >p^{r_s^2} > p^{|b|} >|b|
\end{array}
$$
and $\sigma \not= g$ too. Thus $\mathbf{d}$ is a $T$-sequence.

Let us prove that $s_{\mathbf{d}} (\mathbb{Z} (p) \oplus \mathbb{T} ) = 0\oplus \mathbb{Z} (p^{\infty})$. Let $0\leq \omega <p, x\in \mathbb{T}$ and $(d_n , \omega + x) \to 1$. Then $(d_{2n} , \omega + x)= (p^n ,x)\to 1$. Hence $x\in \mathbb{Z} (p^{\infty})$ (see \cite{Arm} or remark 3.8 \cite{BD2}). Let $x= \frac{\rho}{p^{\tau}} , \rho\in \mathbb{Z}, \tau >0$, then for $n> 2\tau$ we have $(d_{2n-1} , \omega + x)= (e,\omega)\to 1$ only if $\omega =0$.
Hence ${\rm Cl} (s_{\mathbf{d}} (G^{\wedge})) = \mathbb{T}$. By theorem 4 \cite{Ga1},
$\mathbf{n} (G, \mathbf{d}) =\mathbb{Z} (p)$ is finite.

2) $G$ is a torsion group and it contains a subgroup of the form $\mathbb{Z} (p^{\infty})$. Then, by lemma \ref{l3},  we can assume that  $G =\mathbb{Z} (p^{\infty})$. Let $u\in \mathbb{Z} (p^{\infty})$ and $\beta$ its order. G.Luk\'{a}cs \cite{Luk} defined the following sequence $$
d_{2n-1} = u + \frac{1}{p^{n^3 - n^2}} +\dots + \frac{1}{p^{n^3 - 2n}} +  \frac{1}{p^{n^3 - n}} +  \frac{1}{p^{n^3}}, \quad d_{2n} =  \frac{1}{p^{n}}, \quad n\geq 1.
$$
He proved that if $p>2$ then $\mathbf{d} =\{ d_n \}$ is a $T$-sequence. Now we give a simple proof that $\mathbf{d}$ is a $T$-sequence for every $p$.

Let $g=\frac{b}{p^z}\not= 0 \in \mathbb{Z} (p^{\infty}), b\in \mathbb{Z}$, be an irreducible fraction and $k\geq 0$. Then, for some $q\geq z$, $\frac{1}{p^q}$ generates the subgroup containing $u$ and $g$. Set $t= p(k+1)+q$ and $m=10t$. We shall prove that $g\not\in A(k,m)$.

a) Let $\sigma \in A(k,m)$ have the following form
$$
\sigma =l_1 d_{2r_1 } + l_2 d_{2r_2 } +\dots + l_s d_{2r_s} = l_1 p^{r_1}  +\dots + l_s p^{r_s} = p^{r_1} \cdot\sigma' ,
$$
where  $m\leq 2r_1 <2r_2 <\dots <2r_s$. Then
$$
|\sigma | = |l_{1} d_{2r_{1}} +l_{2} d_{2r_{2}} +\dots + l_{h} d_{2r_{h}} | \leq \sum_{i=1}^h \frac{|l_i|}{p^{r_i}} \leq  \frac{k+1}{p^{r_1}} < \frac{k+1}{p^{k+1+q}} < \frac{1}{p^{q}}
$$
and $\sigma \not= g ({\rm mod} 1)$.

b) Let $\sigma \in A(k,m)$ have the following form
$$
\sigma =l_1 d_{2r_1 -1} + l_2 d_{2r_2 -1} +\dots + l_s d_{2r_s -1},
$$
where $m< 2r_1 -1 < 2r_2 -1 < \dots < 2r_s -1 $ and integers $l_1, l_2,\dots , l_s$ be such that $l_s \not= 0$ and  $\sum_{i=1}^s |l_i | \leq k+1$. Since $n^3 < (n+1)^3 -(n+1)^2$ and $r_s >5p(k+1)$, we have
$$
\sigma = \frac{z'}{p^{r_s^3 - r_s}} + \frac{l_s}{p^{r_s^3 }}, \mbox{ where } z' \in \mathbb{Z}.
$$
Since $|l_s | < k+1< \frac{r_s}{p} < p^{r_s -1}$, we have the following: if $\sigma =\frac{z''}{p^{\alpha}}, z'' \in \mathbb{Z}$, is an irreducible fraction, then $\alpha \geq r^3_{s} - r_s +1 > 5q$ and $\sigma \not= g ({\rm mod} 1)$.

c) Let $\sigma \in A(k,m)$ have the following form
$$
\sigma = l_1 d_{2r_1 -1} + l_2 d_{2r_2 -1} +\dots + l_s d_{2r_s -1} +l_{s+1} d_{2r_{s+1}} +l_{s+2} d_{2r_{s+2}} +\dots + l_{h} d_{2r_{h}},
$$
where $0< s<h$ and
$$
\begin{array}{lr}
  m< 2r_1 -1 < 2r_2 -1 < \dots < 2r_s -1 , &   \\
  m\leq 2r_{s+1} < 2r_{s+2} < \dots < 2r_{h}, & l_i \in \mathbb{Z} \setminus \{ 0\} , \; \sum_{i=1}^h |l_i | \leq k+1.
\end{array}
$$
Since the number of summands in $d_{2r_s -1}$ is $r_s +2 > 5p(k+1)+q +2$ and $h-s <k+1$, there exists $r_s - 2> i_0 >2$ such that for every $1\leq w \leq h-s$ we have
$$
\mbox{either } r_{s+w} < r_s^3 - (i_0 +2)r_s  \mbox{ or }  r_{s+w} > r_s^3 - (i_0 -1)r_s .
$$
The set of all $w$ such that $r_{s+w} < r_s^3 - (i_0 +2)r_s$ we denote by $K$ (it may be empty or has the form $\{ 1,\dots ,a\}$ for some $1\leq a \leq h-s$). Set $L =\{ 1,\dots ,h-s \} \setminus B$. Thus
$$
\begin{array}{ll}
  \sigma = & l_1 d_{2r_1 -1}  +\dots + l_{s-1} d_{2r_{s-1} -1} + \sum_{w\in K} l_{s+w} d_{2r_{s+w}} + \frac{l_s}{p^{r^3_{s} - r_s^2 }} +\dots \frac{l_s}{p^{r^3_{s} - (i_0 +2)r_s}} +\\
   & \frac{l_s}{p^{r^3_{s} - (i_0 +1)r_s}} +\frac{l_s}{p^{r^3_{s} - i_0 r_s}}+ \\
   & \frac{l_s}{p^{r^3_{s} - (i_0 -1)r_s}} + \dots +\frac{l_s}{p^{r^3_{s} }} + \sum_{w\in L} l_{s+w} d_{2r_{s+w}} .
\end{array}
$$
Since $n^3 < (n+1)^3 -(n+1)^2$, then the expression in row 1 can be represented in the form $\frac{c}{p^{r^3_{s} - (i_0 +2)r_s}}$, for some $c \in \mathbb{Z}$.
Since $r_s > 5p(k+1)$, then $\frac{1}{1- 1/p^{r^s}} < \frac{32}{31}$ and $2k < p^{2k} <p^{r_s}$. Thus
$$
\left| \left( \frac{l_s}{p^{r^3_{s} - (i_0 -1)r_s}} + \dots +\frac{l_s}{p^{r^3_{s} }} \right) + \sum_{w\in L} l_{s+w} d_{2r_{s+w}} \right| <
\frac{|l_s|}{p^{r^3_{s} - (i_0 -1)r_s}} \cdot \frac{1}{1-\frac{1}{p^{r_s}}} + \frac{k}{p^{p^{r^3_{s} - (i_0 -1)r_s} +1}} <
$$
$$
\frac{1}{p^{r^3_{s} - (i_0 -1)r_s}} \left( k\frac{32}{31} + k\frac{1}{p} \right) <
\frac{2k}{p^{r^3_{s} - (i_0 -1)r_s}} < \frac{1}{p^{r^3_{s} - i_0 r_s}},
$$
and we have the following: if $\sigma =\frac{c''}{p^{\alpha}}, c'' \in \mathbb{Z}$, is an irreducible fraction, then $\alpha \geq r^3_{s} - (i_0 +1)r_s > 5q$ and $\sigma \not= g$. Thus $\{ d_n\}$ is a $T$-sequence.

Now we can repeat the proof of theorem 4.4(b) \cite{Luk}. If $x\in s_{\mathbf{u}} (X)$ and $(x, d_{2n}) \to 1$, then $(x, \chi) = \exp ( 2\pi i m \chi), \forall \chi \in \mathbb{Z} (p^{\infty})$, for some $m\in \mathbb{Z}$ (example 2.6.3 \cite{ZP2}). Since
$$
0 \leq \frac{1}{p^{n^3 - n^2}} +\dots + \frac{1}{p^{n^3 - 2n}} +  \frac{1}{p^{n^3 - n}} +  \frac{1}{p^{n^3}} \leq \frac{n+1}{p^{n^3 - n^2}} \to 0
$$
one has $(x, d_{2n-1}) \to (x,u)$. Thus $(x,u)=1$.

On the other hand, if $x_0 =(0,\dots,0,1,0,\dots)$, where 1 occupies  the position $\beta +1$, then $(x_0 , d_n)\to 1$. Since the closure of $\langle x_0 \rangle$ is $\langle u \rangle^{\perp}$ (remark 10.6 \cite{HR1}), we have $\mathbf{n} (\widehat{X}, \mathbf{u}) = \langle u \rangle$.

3) $G$ is a torsion group and it contains no a subgroup of the form $\mathbb{Z} (p^{\infty})$.
Then, by the Kulikov theorem (corollary 24.3 \cite{Fuc}), there exist some prime $p_1$ and $1\leq k_1 <\infty$ such that
$$
X^{\wedge} =\mathbb{Z} (p_1^{k_1}) \oplus H_1 .
$$
Since $H_1$ is infinite, again by the Kulikov theorem, there exist some prime $p_2$ and $1\leq k_2 <\infty$ such that
$$
H_1 =\mathbb{Z} (p_2^{k_2}) \oplus H_2 .
$$
Continuing this process we can see that $G$ must contain a subgroup of the form $\bigoplus_{i=1}^{\infty} \mathbb{Z} (p_i)$ with prime $p_i$. Thus we can assume that $G$ has this form. By theorem 3.1 \cite{Luk}, $G$ admits an almost maximally-almost periodic group topology.
$\Box$

\begin{rem} \label{r}
I.Protasov and E.Zelenyuk \cite{ZP1} proved that there exists  a $T$-sequence $\mathbf{d}$ such that $(\mathbb{Z},\mathbf{d})$ has only  the trivial character. Now we give a simple {\it concrete } such $T$-sequence. Let $G=\mathbb{Z}$ and $\gamma$ and $q>1$ be any positive integers. Set
$$
d_{2n-1} = \gamma + q^{n^3 - n^2} +\dots + q^{n^3 - 2n} +  q^{n^3 - n} + q^{n^3}, \quad d_{2n} =  q^{n}, \quad n\geq 1.
$$
Then $\mathbf{d} =\{ d_n\}$ is a $T$-sequence. If $\gamma =q$, then $\mathbf{n} (\mathbb{Z},\mathbf{d})= q\mathbb{Z}$. If $\gamma =1$, then $\mathbf{n} (\mathbb{Z},\mathbf{d})= \mathbb{Z}$. (Indeed, exactly the same proof which is in item 1) of the proof of theorem \ref{t1}
(putting $a=0$ and $t= |b| + (k+1)(m+\gamma)$) shows that $\mathbf{d}$ is a $T$-sequence.

Let $x\in \mathbb{T}$ and $(d_n ,x)\to 1$. Then $(d_{2n} ,x)= (q^n ,x)\to 1$. Hence
$$
x\in \mathbb{Z} (t_1^{\infty}) \oplus \dots \oplus \mathbb{Z} (t_{\alpha}^{\infty}),
$$
where $t_1,\dots , t_{\alpha}$ are the primes that divide $q$ (remark 3.8 \cite{BD2}). Hence, for enough large $n$ we have $(d_{2n-1} ,x)= (\gamma ,x)$. Thus: if $\gamma =q$, then $s_{\mathbf{d}} (\mathbb{T}) =\mathbb{Z} (q)$ and $\mathbf{n} (\mathbb{Z},\mathbf{d})= q\mathbb{Z}$; if $\gamma =1$, then $\mathbf{n} (\mathbb{Z},\mathbf{d})= \mathbb{Z}$.) $\Box$
\end{rem}

{\bf Proof of theorem \ref{t2}}. Let $\mathbf{n} (G)$ be dually embedded. Since $\mathbf{n} (G)$ is dually closed, then $\mathbf{n} (\mathbf{n} (G)) = \mathbf{n} (G)$ by  lemma \ref{l3}.
Conversely, if $\mathbf{n} (\mathbf{n} (G)) \not= \mathbf{n} (G)$, then there exists a nonzero character $\chi \in \mathbf{n} (G)^{\wedge}$. If $\chi$ is extended to $\widetilde{\chi} \in G^{\wedge}$, then, by the definition of $\mathbf{n} (G)$, we must have $(\widetilde{\chi} , \mathbf{n} (G)) = (\chi , \mathbf{n} (G)) =1$ and $\chi$ is trivial. It is a contradiction.

Let $G=\mathbb{Z}$. Then any its nontrivial subgroup $H$ has the form $p\mathbb{Z}$ for some positive integer $p$. If $H=\mathbf{n} (\mathbb{Z},\tau)$, then $H$ is closed. Since $G/H$ is finite and discrete, $H$ is open and, hence, dually embedded \cite{Nob}. The assertion follows. $\Box$

{\large Department of Mathematics, Ben-Gurion University of the
Negev,

Beer-Sheva, P.O. 653, Israel}

{\it E-mail address}: $\quad$ saak@math.bgu.ac.il

\end{document}